\documentclass[a4paper,12pt]{amsart}
\usepackage{amssymb,amsthm}

\input xy
\xyoption{all}

\begin{document}

\newcommand{\opp}{\bowtie }
\newcommand{\po}{\text {\rm pos}}
\newcommand{\supp}{\text {\rm supp}}
\newcommand{\End}{\text {\rm End}}
\newcommand{\diag}{\text {\rm diag}}
\newcommand{\Lie}{\text {\rm Lie}}
\newcommand{\Ad}{\text {\rm Ad}}
\newcommand{\car}{\mathcal R}
\newcommand{\Tr}{{\rm Tr }}
\newcommand{\Spec}{\text{\rm Spec}}
\newcommand{\Hom}{\text {\rm Hom}}

\newtheorem*{thm1}{Theorem}
\newtheorem*{thm1-sigma}{Theorem}
\newtheorem*{thm2}{Lemma 3.5}
\newtheorem*{thm4}{Proposition 4.4}
\newtheorem*{thm5}{Corollary 4.5}
\newtheorem*{th5}{Proposition}
\newtheorem*{thm6}{Lemma 5.2}
\newtheorem*{thm-lambda}{Lemma 5.4}
\newtheorem*{thm7}{Corollary 6.3}
\newtheorem*{thm8}{Lemma 7.1}
\newtheorem*{thm9}{Lemma 7.2}
\newtheorem*{thm10}{Theorem 7.3}
\newtheorem*{thm11}{Corollary 7.4}
\newtheorem*{prop}{Proposition 6.2}

\theoremstyle{remark}
\newtheorem*{remark}{Remark}

\def\ge{\geqslant}
\def\le{\leqslant}
\def\a{\alpha}
\def\b{\beta}
\def\c{\chi}
\def\g{\gamma}
\def\G{\Gamma}
\def\d{\delta}
\def\D{\Delta}
\def\L{\Lambda}
\def\e{\epsilon}
\def\et{\eta}
\def\io{\iota}
\def\o{\omega}
\def\p{\pi}
\def\ph{\phi}
\def\ps{\psi}
\def\s{\sigma}
\def\t{\tau}
\def\th{\theta}
\def\k{\kappa}
\def\l{\lambda}
\def\z{\zeta}
\def\v{\vartheta}
\def\x{\xi}
\def\i{^{-1}}

\def\ca{\mathcal A}
\def\cb{\mathcal B}
\def\cc{\mathcal C}
\def\cd{\mathcal D}
\def\ce{\mathcal E}
\def\cf{\mathcal F}
\def\cg{\mathcal G}
\def\ch{\mathcal H}
\def\ci{\mathcal I}
\def\cj{\mathcal J}
\def\ck{\mathcal K}
\def\cl{\mathcal L}
\def\cm{\mathcal M}
\def\cn{\mathcal N}
\def\co{\mathcal O}
\def\cp{\mathcal P}
\def\cq{\mathcal Q}
\def\car{\mathcal R}
\def\cs{\mathcal S}
\def\ct{\mathcal T}
\def\cu{\mathcal U}
\def\cv{\mathcal V}
\def\cw{\mathcal W}
\def\cz{\mathcal Z}
\def\cx{\mathcal X}
\def\cy{\mathcal Y}

\def\tz{\tilde Z}
\def\tl{\tilde L}
\def\tc{\tilde C}

\title[The closure of Steinberg fibers]
{Closures of Steinberg fibers in twisted
wonderful compactifications}
\author{Xuhua He}%
\address{Department of Mathematics, Massachusetts Institute of Technology, Cambridge, MA 02139, USA}%
\email{xuhua@mit.edu}%
\author{Jesper Funch Thomsen}%
\address{Institut for matematiske fag \\ Aarhus Universitet \\ 8000 \AA rhus C, Denmark}
\email{funch@imf.au.dk}%


\begin{abstract}
By a case-free approach we give a precise description
of the closure of a Steinberg fiber within a twisted
wonderful compactification of a simple linear algebraic
group. In the non-twisted case this description was
earlier obtained by the first author.
\end{abstract}
\maketitle

\section{Introduction}

Let $G$ be a simple linear algebraic group over an
algebraically closed field $k$. Let $B$ be a Borel
subgroup of $G$ and let $T \subset B$ denote a
maximal torus. Let $W$ denote the associated Weyl
group and $I$ denote the associated set of simple
roots. For a subset $J$ of $I$ we let $W_J$ denote
the subgroup of $W$ generated by the simple roots
in $J$.

The wonderful compactification $X$ of $G$ (see e.g. \cite{DP},
\cite{Str}), is a smooth projective $(G \times G)$-variety
containing $G$ as an open subset. The $G \times G$-orbits in $X$
are indexed by the subsets $J$ of $I$, and we fix certain base
points $h_J$ for these orbits. Let $\s$ denote a diagram
automorphism of $G$ and let $X_\s$ be the associated twisted
wonderful compactification of $G$, i.e. as a variety $X_\s$ is
just $X$ but the $G \times G$-action is twisted by $\s$ on the
second coordinate. Let $h_{J,\s}$ denote the point in $X_\s$
identified with $h_{\s(J)}$ in $X$. Then the collection
$h_{\s(J)}$, $J \subset I$, are representatives for the $G \times
G$-orbits  in $X_\s$. A {\it $G$-stable piece} in $X_\s$ is then a
locally closed and smooth subvariety in $X$ of the form
$Z^w_{J,\s} = {\rm diag}(G) (B w, 1) h_{J,\s}$, where $w \in
W^{\s(J)}$ is a minimal length coset representative of
$W/W_{\s(J)}$ and $\diag(G)$ denotes the diagonal in $G \times G$.
We then have a decomposition $X_\s=\bigsqcup_{J \subset I, w \in W^J}
Z^w_J$ (see \cite[12.3]{L2} and \cite[1.12]{H2}).

The $G$-stable pieces were first introduced by Lusztig to study
the $G$-orbits and parabolic character sheaves. However, his
original definition was based on some inductive method. The
(equivalent) definition that we used above is due to the first
author in \cite{H1}. What we need in this paper is that the
dimension of $Z^w_ {J,\s}$ is equal to $\dim(G)-l(w)-|I-J|$, where
$l(w)$ is the length of $w$ and $|I-J|$ is the cardinality of the
set $I-J$ (see \cite[8.20]{L2}). More properties about the
$G$-stable pieces can be found in \cite{L2} and \cite{H2}. The
$G$-stable pieces were also used by Evens and Lu in \cite{EL} to
study the Poisson structure and symplectic leaves.

Consider $G$ as a $G \times G$-variety by left and right
translation and define $G_\s$, similar to the definition
of $X_\s$, by twisting the $G$-structure of $G$ on the second
factor by $\s$. A $\s$-conjugacy class in $G_\s$ is then
a $\diag(G)$-orbit in $G_\s$. The set of elements in
$G_\s$ whose semisimple part lies in a fixed $\s$-conjugacy
class is then called a Steinberg fiber of $G_\s$. In
this paper we study the closure of Steinberg fibers
within $X_\s$.

In \cite{L2}, Lusztig gave an explicit description for the closure
of the unipotent variety in the group compactification when
$G=PGL_2$ or $PGL_3$. In \cite{Sp2}, Springer studied the closure
of an arbitrary Steinberg fiber for any connected, simple
algebraic group and obtained some partial results. Based on
their results, the first author obtained an explicit
description of the closure of Steinberg fibers in the
non-twisted case. The result in \cite{H1} was formulated
using  $G$-stable pieces and the proof was based on a
case-by-case checking.  The main purpose of this paper
is to generalize the result of \cite{H1} to the
twisted case with a more conceptual (and easier) proof.
More precisely, we prove,

\begin{thm1-sigma} Let $F$ be a Steinberg fiber
of $G_\sigma$ and $\bar{F}$ be its closure in
$X_\sigma$. Then
$$\bar{F}-F=\bigsqcup_{J \subset I} \bigsqcup_{\substack{w \in W^{\sigma(J)}
\\  \supp_\sigma(w)=I}}
Z^w_ {J,\sigma},$$
where $\supp_\sigma(w)$ denotes the minimal $\s$-stable
subset of $I$ such that $w$ is contained in
$W_{\supp_\sigma(w)}$.
\end{thm1-sigma}

As a consequence, the boundary of the closure is independent of
the choice of the Steinberg fiber. Likewise it may be shown
that the boundary of the closure of $F$ within any
equivariant embedding of $G$ is independent of the choice of
$F$ (see \cite{T}). As a by-product, we will also give an explicit
description of the ``nilpotent cones'' on $X$.

\section{Wonderful compactifications and $G$-stable pieces}

\subsection*{2.1} Let $G$ denote a simple linear
algebraic group over an algebraically closed field
$k$. We consider $G$ as a $G \times G$-variety by left and
right translation. Let $B$ be a Borel subgroup of $G$, $B^-$ be an opposite
Borel subgroup and $T=B \cap B^-$. The unipotent radical
of $B$ (resp. $B^-$) will be denoted by $U$ (resp. $U^-$).
Let $R$ denote the set of roots defined by $T$ and let
$R^+$ denote the set of positive roots defined by $B$.
Let $(\a_i)_{i \in I}$ be the set of simple roots.
For $i \in I$, we denote by $\o_i$ and $s_i$ the fundamental
weight and the simple reflection corresponding to $\a_i$.

We denote by $W$ the Weyl group associated
to $T$. For any subset $J$ of $I$, let $W_J$ be the subgroup of $W$
generated by $\{s_j \mid j \in J\}$ and $W^J$ be the set of minimal
length coset representatives of $W/W_J$.

For $J \subset I$, let $P_J \supset B$ be the standard parabolic
subgroup defined by $J$ and let $P^-_J \supset B^-$ be the parabolic
subgroup opposite to
$P_J$. Set $L_J=P_J \cap P^-_J$. Then $L_J$ is a Levi subgroup of
$P_J$ and $P^-_J$. The semisimple quotient of $L_J$ of adjoint type
will be denoted by $G_J$. We denote by $\p_{P_J}$ (resp.
$\p_{P^-_J}$) the projection of $P_J$ (resp. $P^-_J$) onto $G_J$.

\subsection*{2.2}
Let $X$ denote the wonderful compactification of $G$
(\cite{DP}, \cite{Str}). Then $X$ is an
irreducible, smooth projective $(G \times G)$-variety with
finitely many $G \times G$-orbits $Z_J$ indexed by the subsets $J$
of $I$. As a $(G \times G)$-variety the orbit $Z_J$ is uniquely
isomorphic to the product $(G \times G) \times_{P^-_J \times P_J} G_J$, where
$P^-_J \times P_J$ acts on $G \times G$ by $(q, p) \cdot (g_1,
g_2)=(g_1 q \i, g_2 p \i)$ and on $G_J$ by $(q, p)
\cdot z=\p_{P^-_J}(q) z \p_{P_J}(p) \i$. Let $h_J$ be the image of
$(1, 1, 1)$ in $Z_J$ under this isomorphism.

We denote by ${\rm diag}(G)$ the image of the diagonal embedding
of $G$ in $G \times G$. For $J \subset I$ and $w \in W^J$, set
$Z^w_J={\rm diag}(G) (B w, 1) h_J$. Then $Z^w_J$ is a locally
closed subvariety of $X$ and (see \cite[12.3]{L2} and
\cite[1.12]{H2})
$$X=\bigsqcup_{\substack{J \subset I \\ w \in W^J}} Z^w_J.$$
We call $Z^w_J$ a $G$-stable piece.

\section{Twisted actions}

\subsection*{3.1}
An automorphism $\s$ of $G$ which stabilizes the
Borel subgroup $B$ and the maximal torus $T$ will
induce a permutation of $I$. When the order of
$\s$ as an automorphism of $G$ coincides with the
order of the associated permutation of $I$, we
say that $\s$ is a diagram automorphism. From now
on $\s$ will denote a diagram automorphism of $G$.
We also denote by $\s$ the corresponding bijection
on $I$ and $W$.

Let $G_{\s}$ be the $(G \times G)$-variety which as a variety is
isomorphic to $G$ and where the $G \times G$ action is twisted by
the morphism $G \times G \rightarrow G \times G$, $(g, h) \mapsto
(g, \s(h))$ for $g, h \in G$. Then we define the wonderful
compactification $X_{\s}$ of $G_{\s}$ to be the $G \times
G$-variety which as a variety is isomorphic to the wonderful
compactification $X$ of $G$ and where the $G \times G$ action is
twisted in the same way as above. Notice that we may regard
$G_{\s}$ as a connected component of the semidirect product $G
\rtimes <\s>$. In this case, $X_{\s}$ is the completion of
$G_{\s}$ considered in \cite[12]{L2}.

The $G \times G$-orbits in ${X}_{\s}$ coincide
with the associated orbits in ${X}$ and we let $Z_{J, \s}$ denote
the orbit coinciding with $Z_{\s(J)}$. Accordingly we let $h_{J,
\s}$ denote the point in  $Z_{J, \s}$ identified with the base
point $h_{\s(J)}$ of $Z_{\s(J)}$. For $J \subset I$ and $w \in W^{\s(J)}$, set
$Z^w_{J, \s}={\rm diag}(G) (B w, 1) h_{J, \s}$. Then
$$X_{\s}=\bigsqcup_{J \subset I} \bigsqcup_{w \in W^{\s(J)}} Z^w_{J,
\s}.$$ We call $(Z^w_{J, \s})_{J \subset I, w \in W^{\s(J)}}$ the
$G$-stable pieces of $X_{\s}$ (see \cite[12.3]{L2} and
\cite[1.12]{H2}).

\subsection*{3.2} The orbits of ${\rm diag}(G)$ on $G_{\s}$
are called $\s$-conjugacy classes. Let $G //_{\s} G$ be the
affine variety whose algebra is the subalgebra $k[G]^{G, \s}$
of functions in $k[G]$ invariant under $\s$-conjugacy. The
inclusion $k[G]^{G, \s} \rightarrow k[G]$
induces a morphism ${\rm St}: G_{\s} \rightarrow G//_{\s} G$. If
$\s$ is trivial, then ${\rm St}$ is just the Steinberg morphism of
$G$. Thus for arbitrary $\s$, we call ${\rm St}$ the Steinberg
morphism of $G_{\s}$ and the fibers the Steinberg fibers of
$G_{\s}$.

An element $g \in G_{\s}$ is $\s$-conjugate to an element in $B$
\cite[Lem.7.3]{Ste}.
Write $b = t u$ where $t \in T$ and $u$ is an element of the
unipotent radical $U$ of $B$. It is then easily seen that there
exists an element $t_1 \in T$, such that $t_1 t \s(t_1) \i \in
T^{\s}$. Hence, $g$ is $\s$-conjugate to some element in $T^{\s}
U$, i.e. we may assume that $t \in T^\s$. Notice, that $t$ is
contained in the closure of the $\s$-conjugacy class of $t u$
and thus, by geometric invariant theory, we find ${\rm St}(t u)
= {\rm St}(t)$.
Moreover, considering $t \sigma$ as an element of the semisimple
group  $G \rtimes <\s>$ it follows that $t \sigma$ is
quasi-semisimple in the sense of \cite[Sect.9]{Ste}, i.e.
the automorphism of $G$ obtained by conjugation by $t \sigma$
will fix a Borel subgroup and a maximal torus thereof.
As a consequence, the $\s$-conjugacy class of
$t$ in $G_{\s}$ is closed \cite[II.1.15(f)]{Spa}. We
conclude that any Steinberg fiber of $G_{\s}$ is of the form
$\bigcup_{g \in G} g (t U) \s(g) \i$ for some $t \in T^{\s}$. In
particular, any Steinberg fiber is irreducible.

\subsection*{3.3} Let $G_{\rm sc}$ be the connected,
simply connected, linear algebraic group associated to
$G$, and let $B_{\rm sc}$ (resp. $T_{\rm sc}$) denote the
Borel subgroup (resp. maximal torus) of $G_{\rm sc}$
associated to $B$ (resp. $T$). By \cite[9.16]{Ste} the
automorphism $\s$ of $G$ may be lifted to an automorphism
of $G_{\rm sc}$, which we also denote by $\sigma$.
We then define the $G_{\rm sc} \times G_{\rm sc}$-variety
$G_{{\rm sc},\s}$ similar to the definition of $G_\s$.
We may also form the quotient $G_{\rm sc} //_\s G_{\rm sc}$
and define Steinberg fibers in $G_{{\rm sc},\s}$ similar
to the considerations in 3.2 for $G_\s$.

The automorphism $\s$ of $G_{\rm sc}$ induces a
natural action of $\s$ on the set $\L$ of
$T_{\rm sc}$-characters, and we let $\L^{\s}_+$
denote the set of $\sigma$-invariant dominant weights.
Let $\cc_1, \cc_2, \cdots, \cc_l$ denote the
$\s$-orbits in $I$, and set  $\o_{\cc_j}=
\sum_{i \in \cc_j} \o_i$. Then the elements
$\o_{\cc_j}$, $j=1, \dots,l$,
is a generating set of the semigroup  $\L^{\s}_+$.

\subsection*{3.4}

To any dominant weight $\l =\sum_{i \in I}
a_i \o_i$ we associate the subset $I(\l)=\{i \in I
\mid a_i \neq 0\}$ of $I$. For $w \in W$,
let $\supp(w) \subset I$ be the set of simple roots
whose associated simple reflections occur in
some (or equivalently, any) reduced
decomposition of $w$ (see \cite[Prop.IV.1.7]
{Bou}) and let $\supp_{\s}(w)=\bigcup_{k \ge 0}
\s^k(\supp(w))$. We have the following
characterization of $\supp(w)$.

\begin{thm2} Let $w \in W$ and $i \in I$. Then $w \o_i \neq \o_i$
if and only if $i \in supp(w)$. Hence for a dominant weight $\l$, $w
\l \neq \l$ if and only if $I(\l) \cap \supp(w) \neq \varnothing$.
\end{thm2}

Proof. If $i \notin \supp(w)$, then $w \o_i=\o_i$. Now we fix a
reduced expression $w=s_{i_1} \cdots s_{i_n}$. Assume that $i \in
{\rm supp}(w)$. We show that $w \o_i \neq \o_i$ by induction on $n$. If
$i_n \neq i$, then we are done by induction in $n$. Hence, we may
assume that $i_n=i$. But then $w \alpha_i$ is a negative root.
Thus $1 = \langle \omega_i, \alpha_i^\vee \rangle = \langle w
\omega_i, (w \alpha_i)^\vee \rangle $ and, in particular, we
cannot have $ w \o_i = \o_i$. \qed

\subsection*{3.6}
For any dominant weight $\l \in \L_{+}$ let ${\rm H}({\l})$ denote the
dual Weyl module for $G_{\rm sc}$ with lowest weight $-\l$. We then
define $^\s{{\rm H}({\l})}$ to be the $G_{\rm sc}$-module which as a vector
space is ${{\rm H}({\l})}$ and with $G_{\rm sc}$-action twisted by the
automorphism $\s$ of $G_{\rm sc}$. Notice that up to a nonzero
constant there exists a unique $G_{\rm sc}$-isomorphism $^\s{{\rm H}({\l})}
\simeq {\rm H}(\s(\l))$. In particular, when $\l \in \L^{\s}_+$ is
$\s$-invariant there exists a $G_{\rm sc}$-equivariant isomorphism $f_\l :
{{\rm H}({\l})} \rightarrow ~ ^\s {\rm H}(\l)$. Fix $f_\l$ such that
the its restriction to the lowest weight space $k_{-\lambda}$
in ${{\rm H}({\l})}$ is the identity map (here we use the
identification of $^\s {\rm H}(\l)$ with  ${{\rm H}({\l})}$
as vector spaces).

\section{The ``nilpotent cone'' of $X$}

\subsection*{4.1}
For any dominant weight $\l$ there exists (see \cite[3.9]{DS}) a
$G \times G$-equivariant morphism
$${\rho_{\l}} : X  \rightarrow \mathbb P \bigl( \End({\rm H}({\l})) \bigr)$$
which extends the morphism  $G \rightarrow \mathbb P \bigl(
\End({\rm H}({\l})) \bigr)$ defined by $g \mapsto g [ {\rm
Id}_{\l}]$, where $ [{\rm Id}_{\l}]$ denotes the class
representing the identity map on ${\rm H}({\l})$ and $g$ acts by
the left action. By the definition of $X_{\s}$ we obtain a $G
\times G$-equivariant morphism
$$  X_\s  \rightarrow \mathbb P \bigl( \Hom_k(^\s {\rm H}({\l}),
{{\rm H}({\l})} ) \bigr).$$
When $\l \in \L^{\s}_+$ we may apply $f_\l$ to
obtain an induced map
$$\rho_{\l,\s} : X_\s  \rightarrow
\mathbb P \bigl( \End({\rm H}({\l}) ) \bigr)$$
which is $G \times G$-equivariant.

\subsection*{4.2}
An element in $\mathbb{P} \bigl( \End({\rm H}({\l}) ) \bigr)$ is
said to be nilpotent if it may be represented by a nilpotent
endomorphism of ${\rm H}(\l)$. For  $\l \in  \L^{\s}_+$
we let
$$\cn(\l)_{\s}=\{z \in X_{\s} \mid {\rho_{\l,\s}}(z)
\hbox{ is nilpotent} \},$$
and call $\cn(\l)_{\s}$ the nilpotent cone of $X_{\s}$ associated
to the dominant weight $\l$. In 4.4, we will give an explicit
description of $\cn(\l)_{\s}$.

\subsection*{4.3} Define ${\rm ht}$ to be the height map on
the root lattice, i.e., the linear map on the root lattice
which maps all the simple roots to $1$.

Now fix $\l \in \L_+$. Choose a basis $v_1,...,v_m$ for
${\rm H}({\l})$ consisting of $T$-eigenvectors with eigenvalues
$\lambda_1,...,\lambda_m$ satisfying ${\rm ht}( \l_j + \l)
\geq {\rm ht}(\l_i + \l)$ whenever $j \leq i$. Then $B$ is upper
triangular with respect to this basis.

Let $A_J$ be a representative of ${\rho}_{\l}(h_{J})$ in
$\End({\rm H}({\l}))$. Then when  $\l_j+\l$ is a linear
combination of the simple roots in $J$ we have that $A_J v_j \in
k^\times v_j$. If $\l_j+\l$ is not a  linear combination of the
simple roots in J then $A_J v_j = 0$. Assuming that $\l$ is
$\s$-invariant we obtain, by the definitions in 4.1, a similar
description for a representative $A_{J,\s}$ of
${\rho}_{\l,\s}(h_{J,\s})$ : if $\l_j+\l$ is a linear combination
of the simple roots in J then we have that $A_{J,\s} v_j  \in
k^\times f_\l(v_j)$; otherwise $A_{J,\s} v_j = 0$. Notice that we
regard $f_\l(v_j)$ as an element of ${\rm H}(\l)$ and, as such,
$f_\l(v_j)$ is a $T$-eigenvector of weight $\s(\l_j)$.

\

We now obtain.

\begin{thm4} Let $\l \in \L^{\s}_+$, then $$\cn(\l)_{\s}=\bigsqcup_{J
\subset I} \bigsqcup_{\substack{w \in W^{\s(J)} \\
I(\l) \cap \supp(w) \neq \varnothing}} Z^w_{J, \s}.$$
\end{thm4}

Proof. Let $w \in W^{\s(J)}$. Assume that $w \l \neq \l$. Note
that if $x$ is a nonnegative linear combination of the simple roots
in $J$ then ${\rm ht}(w \s(x)) \ge
{\rm ht}(x)$. Hence,
$${\rm ht} \bigl(w \s(-\l+x)+\l
\bigr)={\rm ht}(w \s(x))+{\rm ht}(-w \l+\l)>{\rm ht}(x).$$
Therefore, $\rho_{\l,\s}((w, 1) h_{J,
\s})$ is represented by a strictly upper triangular matrix with
respect to the chosen basis in 4.3 above. As a consequence for any
$b \in B$, $\rho_{\l,\s}((b w, 1) h_{J,
\s})$ is also represented by a strictly upper triangular matrix. So
$(B w, 1) h_{J, \s} \subset \cn(\l)_{\s}$. Since $\cn(\l)_{\s}$ is
${\rm diag}(G)$-stable it follows
$Z^w_{J, \s}={\rm diag}(G) (B w, 1) h_{J, \s} \subset
\cn(\l)_{\s}$.

Now assume that $w \l=\l$. Let $b \in B$ and $z=(b w, 1) h_{J,
\s}$. Denote by $A$ a representative of ${\rho}_{\l, \s}(z)$ in
$\End({\rm H}({\l}))$. Let $V$ be the subspace of ${\rm H}(\l)$
spanned by $v_1, \cdots, v_{m-1}$. Then $A v_m \in k^{\times}
v_m+V$ and $A V \subset V$. Hence, $A^n v_m \neq 0$ for all $n \in
\mathbb N$. Thus $z \notin \cn(\l)_{\s}$. \qed

\begin{thm5} Let $\l, \mu \in \L^{\s}_+$, then
$$\cn(\l+\mu)_{\s}=\cn(\l)_{\s} \cup \cn(\mu)_{\s}.$$
\end{thm5}

Proof. This follows from the relation $I(\l+\mu)=I(\l) \cup
I(\mu)$. \qed

\section{A compactification of $G_{\rm sc}$}

\subsection*{5.1}
Consider the morphism $\psi_i: G_{\rm sc} \rightarrow {\mathbb P}
\bigl(\End({\rm H}(\o_{i})) \oplus k \bigr)$ defined by $g \mapsto
[(g \cdot {\rm Id}_{{\rm H}(\o_{i})},1)]$, where ${\rm Id}_{{\rm H}(\o_{i})}$
denotes the identity map on ${\rm H}(\o_{i})$ and $g$ acts on
$\End({\rm H}(\o_{i}))$ by the left action. Let furthermore
$\pi: G_{\rm sc} \rightarrow X$ denote the the natural $G_{\rm
sc} \times G_{\rm sc}$-equivariant morphism and let $X_{\rm sc}$
denote the closure of the image of the product map
$$(\pi, \prod_{i \in I} \psi_i): G_{\rm sc}
\rightarrow X \times \prod_{i \in I} {\mathbb P}
\bigl(\End({\rm H}(\o_{i}) \oplus k),
$$
Then $X_{\rm sc}$ is a projective
$G_{\rm sc} \times G_{\rm sc}$-equivariant
variety containing $G_{\rm sc}$ as an open subset. Unlike $X$ the
variety $X_{\rm sc}$ need not be smooth and in general it is not
even normal (e.g. in type $A_3$). By abuse of
notation we use the notation $\pi$ and $\psi_i$,
$i \in I$, for the natural extensions of
the corresponding maps to $X_{\rm sc}$.

\begin{thm6} The projective morphism $\pi: X_{\rm sc}
\rightarrow X$ defines a bijection between
$X_{\rm sc}-G_{\rm sc}$ and $X-G$. In particular,
$\pi$ is a finite morphism.
\end{thm6}

Proof. Notice that the $G_{\rm sc} \times G_{\rm sc}$-invariant
homogeneous polynomial function on $\End({\rm H}(\o_{i}))
\oplus k$ defined by $(f,a) \mapsto {\rm det}(f) -
a^{{\rm dim}_k({\rm H}(\o_{i} ))}$ vanishes at
$({\rm Id}_{{\rm H}(\o_{i})},1)$. Thus
$$X_{\rm sc} \subset X \times \prod_{i \in I}  \biggr({\mathbb P}
\bigl(\End({\rm H}(\o_{i}) \oplus k) - {\mathbb P} (0 \oplus k)
\biggr),$$
and we may consider the following commutative diagram
$$\xymatrix{ X_{\rm sc} \ar[d]^{\pi} \ar[r]^(.15){} &
X \times \prod_{i \in I}  \biggr({\mathbb P}
\bigl(\End({\rm H}(\o_{i}) \oplus k)
- {\mathbb P} (0 \oplus k) \biggr) \ar[d] \\
X \ar[r]^(.3){{\rm id}_X \times \prod_{i \in I}
{\rho}_{\o_{i}} }
& X \times
\prod_{i \in I}  {\mathbb P}
\bigl(\End({\rm H}(\o_{i})) \bigr)   \\
},
$$
where all the maps are the natural ones. Assume now that $x$ is an
element of the boundary $X_{\rm sc}-G_{\rm sc}$. As the dimensions
of $G_{\rm sc}$ and $X_{\rm sc}$ coincide the $(G_{\rm sc}, 1)$-stabilizer
of $x$ has strictly positive dimension. In particular, the images
$\psi_i(x) = [(f_i,a_i)]$, $i \in I$,  have the same property.
Thus, the endomorphism $f_i$ is not invertible and thus $a_i=0$.
This proves that
$$X_{\rm sc} - G_{\rm sc} \subset
X \times \prod_{i \in I}  {\mathbb P} \bigl(\End({\rm H}(\o_{i}) ),
$$ and hence $\pi$ maps $X_{\rm sc}-G_{\rm sc}$ injectively to the
boundary $X - G$. As $\pi$ is dominant and projective, and thus
surjective, this proves the first assertion. Finally $\pi$
is finite as it is projective and quasifinite.
\qed

\subsection*{5.3} For a dominant weight
$\l \in \L_+$ let
$\psi_\l : G_{\rm sc} \rightarrow {\mathbb P}
\bigl(\End({\rm H}(\l)) \oplus k \bigr)$
be the morphisms defined by $\psi_\l(g) =
[(g \cdot {\rm Id}_{{\rm H}(\l)},1)]$,
for $g \in G_{\rm sc}$. Then we let
$X_{\rm sc}^\lambda$ denote the closure
of the image of the map
$$ G_{\rm sc} \rightarrow X_{\rm sc} \times
{\mathbb P} \bigl(\End({\rm H}(\l)) \oplus k \bigr)$$
defined as the product of the inclusion
$G_{\rm sc} \subset X_{\rm sc}$ and $\psi_\l$.
We claim

\begin{thm-lambda}
The canonical morphism $\pi^\l : X_{\rm sc}^\l
\rightarrow X_{\rm sc}$ is an isomorphism.
\end{thm-lambda}
Proof. Let $X_0$ be the complement of the union
of the closures $\overline{B s_i B^-}$, $i \in I$,
within $X$, and let $X_0' = \overline{T} \cap X_0$.
Then, by \cite[Prop.6.2.3(i)]{BK}, the natural
morphism
$$ U \times U^- \times X_0' \rightarrow X_0,$$
$$ (g,h,x) \mapsto (g,h)\cdot x,$$
is an isomorphism of varieties. Let
$X_{{\rm sc},0}=\pi^{-1}(X_0)$ and
$X_{{\rm sc},0}'$ denote the (scheme theoretic)
inverse image $ \pi^{-1}(X'_0)$. As $\pi$ is
$G_{\rm sc} \times G_{\rm sc}$-equivariant we obtain an induced
isomorphism
$$  U \times U^- \times X_{{\rm sc},0}'
\rightarrow X_{{\rm sc},0}.$$
In particular, $X_{{\rm sc},0}'$ is an
irreducible closed subvariety of the open
subvariety $X_{{\rm sc},0}$ containing
$T_{\rm sc}$ as an open subset. Thus,
$X_{{\rm sc},0}'$ is contained in the
closure of $T_{\rm sc}$ within $X_{\rm sc}$.
Let $\pi_\l$ denote the composition of $\pi$
and $\pi^\l$. Defining $X^\l_{{\rm sc},0}=\pi_\l^{-1}
(X_0)$ and $(X^\l_{{\rm sc},0})'= \pi_\l^{-1}(X'_0)$
we, similarly, obtain an isomorphism
$$  U \times U^- \times (X_{{\rm sc},0}^\l)'
\rightarrow X^\l_{{\rm sc},0}.$$
Notice that $X^\l_{{\rm sc},0} =
\pi^\l( X_{{\rm sc},0})$. Moreover,
the $G \times G$-translates of $X_0$ cover
$X$ \cite[Thm.6.1.8]{BK}. Thus, it
suffices to show that the morphism
$(X_{{\rm sc},0}^\l)' \rightarrow
X_{{\rm sc},0}'$ induced by $\pi^\l$
is an isomorphism. This will follow
if $\pi^\l$ induces an isomorphism
between the closures of $T_{\rm sc}$
in $X_{\rm sc}$ and $X_{{\rm sc}}^\l$.
Determining the latter closures of
$T_{\rm sc}$ and checking that they
are isomorphic is now an easy exercise.
\qed

\

It follows that we may consider $\psi_\l$
as the extended morphism
$$\psi_\l : X_{\rm sc} \rightarrow {\mathbb P}
\bigl(\End({\rm H}(\l)) \oplus k \bigr),$$
which we will do in the following. As in
the proof of Lemma 5.2 we may prove
that
$$\psi_\l(X_{\rm sc}) \subset
\biggl({\mathbb P} \bigl(\End({\rm H}(\l)) \oplus k \bigr) -
{\mathbb P}\bigl(0 \oplus k \bigr)\biggr)$$ and that the induced
map $X_{\rm sc} \rightarrow {\mathbb P}\bigl(\End({\rm H}(\l))
\bigr)$ is compatible with $\pi : X_{\rm sc} \rightarrow X$ and
the map $\rho_\l : X \rightarrow {\mathbb P}\bigl(\End({\rm
H}(\l)) \bigr)$.

\subsection*{5.5} The variety $X_{\rm sc}$ is a compactification of
$G_{\rm sc}$ with the $G_{\rm sc} \times G_{\rm sc}$ action defined in the
natural way. Let $X_{\rm sc, \s}$ be the $G_{\rm sc} \times
G_{\rm sc}$-variety which as a variety is isomorphic to $X_{\rm sc}$ and
where the $G_{\rm sc} \times G_{\rm sc}$-action is twisted by the morphism
$G_{\rm sc} \times G_{\rm sc} \rightarrow G_{\rm sc} \times G_{\rm sc}$ ,
$(g,h) \mapsto (g, \s(g))$ for $g,h \in G_{\rm sc}$. Thus
we may identify $G_{{\rm sc},\s}$ of 3.3 with an open
subset of $X_{{\rm sc},\s}$. Notice
that by Lemma 5.2
we may identify the boundaries of $X_{\rm sc,\s}$
and $X_\s$ and we may therefore also regard $Z_{J,\s}^w$, for
$J \neq I$, as subsets of $X_{\rm sc,\s}$.

\section{Steinberg fibers and trace maps}

\subsection*{6.1} Let $\Tr_i$ denote the
trace function on $\End({\rm H}({\o_{\cc_i}}))$. To each $a_i \in
k$ we may associate a global section $(\Tr_i, a_i)$ of the line
bundle $\co_i(1):=\co_{{\mathbb P} \bigl(\End({\rm H}(\o_{\cc_i}))
\oplus k \bigr)}(1)$ on ${\mathbb P} \bigl(\End({\rm
H}(\o_{\cc_i})) \oplus k \bigr)$. The pull back of $(\Tr_i, a_i)$
to $X_{\rm sc, \s}$, by the morphism $\psi_{\o_{\cc_i}}$, is then
a global section $f^{\s}_{i, a_i}$ of a line bundle on $X_{\rm sc,
\s}$. In the following, we will study the common zero set $Z(a_1,
\cdots, a_l)$ of the sections $f^{\s}_{i,a_i}$, for varying $a_i
\in k$. By choosing a trivialization of the pull back of
$\co_i(1)$ to $G_{\rm sc, \s}$ we may think of $f^{\s}_{i, a_i}$
as a function on $G_{\rm sc, \s}$ and, by abuse of notation, we
also denote this function by $f^{\s}_{i, a_i}$. We assume that the
trivialization is chosen such that $f^{\s}_{i, a}- f^{\s}_{i, 0}=
a$ as functions on $G_{\rm sc}$ (then the trivialization is
actually unique). Then $f^{\s}_{i, a_i}$ is invariant under
$\s$-conjugation by $G_{\rm sc}$ and thus $f^{\s}_{i, a_i}$
induces a morphism $\bar{f}^{\s}_{i,a_i} : G_{\rm sc} //_{\s}
G_{\rm sc} \rightarrow k$. We then claim

\begin{prop}
The product morphism
$$(\bar{f}^{\s}_{1,0},\bar{f}^{\s}_{2,0},
\dots,  \bar{f}^{\s}_{l,0}) :
G_{\rm sc}//_\s G_{\rm sc} \rightarrow \mathbb{A}^l,$$
is an isomorphism.
\end{prop}
Proof. Let $f^{\s}_{i}$ denote the restriction
of $f^{\s}_{i, 0}$ to $T_{\rm sc}$. An easy
calculation shows that $f^{\s}_{i}$ equals
$$-w_0 \cc_i + \sum_{\substack{\lambda \in \L^\s, \l \neq \cc_i \\
{\rm H}(\cc_i)_\l \neq 0}} q_{i,\l} \l.$$
Thus  $f^{\s}_{i}$ is contained in the semigroup
algebra $k[\L^\sigma]$ generated by the $\s$-invariant
weights of $\L$. Moreover, $f^{\s}_{i}$
is invariant under the action of the group $W^\s$ of
$\s$-invariant elements of $W$. Hence,
$f^{\s}_{i}$ is an element of the polynomial
ring $k[\L^\sigma]^{W^\s}$ in $l$ variables
\cite[proof of Cor.2]{Sp3}, generated
(as an $k$-algebra) by the elements
$$ {\rm sym}(\cc_i) := \sum_{w \in W^\sigma}
w \cc_i, ~i=1, \dots, l.$$
As in \cite[proof of Lem.6.3]{Ste2} we
conclude that also $f^\s_i$, $i=1, \dots,l$,
generates $k[\L^\sigma]^{W^\s}$
as an $k$-algebra. Now apply \cite[Thm.1]{Sp3}.
\qed

\begin{thm7}
The intersection of $Z(a_1, \cdots, a_l)$ with the boundary
$X_{\rm sc, \s} - G_{\rm sc, \s}$ of $X_{\rm sc, \s}$ is
independent of $a_1,\dots,a_l$. Moreover, the intersection $Z(a_1,
\cdots, a_l) \cap G_{\rm sc, \s}$ is a single Steinberg fiber.
\end{thm7}

Proof. As in the proof of Lemma 5.2 it follows that $x$ is an
element of $X_{\rm sc, \s} - G_{\rm sc, \s}$ exactly when the
image $\psi_{\o_{\cc_i}}(x)$ is of the form $[(f,0)]$ for all
$i=1, \dots,l$. Thus, the section ${f}^{\s}_{i,a_i}$ coincides
with ${f}^{\s}_{i,0}$ on the boundary of $X_{\rm sc, \s}$. This
proves the first statement. The latter statement follows by
Proposition 6.2. \qed

\section{Proofs of the main results}

\begin{thm8} Let $J \subsetneq I$, $w \in W^{\s(J)}$ and $b \in
B$. If ${f}^{\s}_{i,0}((b w,1) h_{J, \s})=0$, then
either (1) $w \o_{\cc_i} \neq \o_{\cc_i}$ or (2) $\cc_i \subset J$
and $w \a_j= \a_j$ for all $j \in \cc_i$.
\end{thm8}

Proof. Assume that $w \o_{\cc_i}=\o_{\cc_i}$. Then the diagonal
entry of a representative $A$ of ${\rho}_{\o_{\cc_i},\s}((b w,1)
h_{J, \s})$ associated to the lowest weight space is nonzero. In
particular, the relation ${f}^{\s}_{i,0}((b w,1) h_{J, \s})=0$
cannot be satisfied unless there exists a weight $x-\o_{\cc_i} \neq
\o_{\cc_i}$ of ${\rm H}({\o_{\cc_i}})$ satisfying $ x=\sum_{j \in J}
a_j \a_j$, with $a_j \in \mathbb N \cup \{0\}$, and $w \s(x)=x$.

Let $K \subseteq J$ denote the set of $j \in J$ such that $a_j
\neq 0$. As $x-\o_{\cc_i}$ is a weight of ${\rm H}({\o_{\cc_i}})$
we know that $\cc_i \cap K$ is nonempty. Now $\sum_{j \in K} a_j w
\a_{\s(j)} =\sum_{j \in K} a_j \a_j$ and thus $\sum_{j \in K} a_j
({\rm ht}(w \a_{\s(j)})-{\rm ht}(\a_j))=0$. As $w \in W^{\s(J)}$
we conclude that  ${\rm ht}(w \a_{\s(j)}) \ge 1$ and consequently
$w \a_{\s(j)}$ is a simple root for all $j \in K$. By the
assumption $w \o_{\cc_i}=\o_{\cc_i}$ we know that $w \a_{\s(j)}=
\a_{\s(j)}$ for each $j \in \cc_i \cap K$. In particular, when $j
\in \cc_i \cap K$ then $a_{\s(j)} = a_j$. Hence, $\cc_i \cap K$ is
invariant under $\s$ and as $\cc_i$ is a single $\s$-orbit
we have $\cc_i \cap K = \cc_i$. This ends the proof. \qed

\begin{thm9} Let $J \subsetneq I$.
Then
$$Z(a_1, \cdots, a_l) \cap Z_{J, \s}= \bigsqcup_{\substack{
w \in W^{\s(J)} \\ \supp_{\s}(w)=I}} Z^w_{J, \s}.$$
\end{thm9}

Proof. By Corollary 6.3 it suffices to consider the case when all
$a_i$ are zero. By Proposition 4.4,
$$\bigsqcup_{J \subset I} \bigsqcup_{\substack{w \in W^{\s(J)} \\
\supp_{\s}(w)=I}}  Z^w_{J, \s}=\cap_{i}
\cn(\o_{\cc_i})_{\s} \subset Z(0, \cdots, 0).$$
This proves one inclusion. For $z \in Z(0,
\cdots, 0) \cap Z_{J, \s}$, we have that $z=(g, g)(b w,1)
h_{J,\s}$ for some $g \in G, b \in B$ and $w \in W^{\s(J)}$. Then
${f}^{\s}_{i,0}((b w,1) h_{J, \s})=0$ for all $i=1, \dots,l$. It suffices to
prove that $\supp_{\s}(w)=I$.

If $w=1$, then by Lemma 7.1, $\cc_i \subset J$ for each $\s$-orbit
$\cc_i$. Thus $I=J$, which contradicts our assumption. Now assume
that $w \neq 1$ and that $\supp_{\s}(w) \neq I$. As $G$ is assumed
to be simple, there exist
simple roots $\a_i$ and $\a_j$ with $n=-\langle \a_{j},\a_{i}^\vee
\rangle \neq 0$ satisfying that $i \in \supp_{\s}(w)$ and $j
\notin \supp_{\s}(w)$. Let $\cc_i$ and $\cc_j$ denote the
associated $\s$-orbits of $\a_i$ and $\a_j$. As
$\supp_{\s}(w)$ is  $\s$-stable it follows that $\cc_i \subset
\supp_{\s}(w)$ and $\cc_j \subset I - \supp_{\s}(w)$.

Now there exists $m \in \mathbb N$, such that $\s^m(i) \in
\supp(w)$ and thus $w \o_{\s^m(i)} \neq \o_{\s^m(i)}$. Hence
redefining, if necessary, $\a_i$ and $\a_j$, we may assume that $w
\o_i \neq \o_i$. Consider then the relation $\a_j = 2 \o_j - n
\o_i - \l$ with $\l$ denoting a dominant weight. Then Lemma 7.1,
applied to $\cc_j$,
implies that $w \a_j = \a_j$ and $w \o_j = \o_j$ and thus $w (n
\o_i + \l) = n \o_i + \l$. As both $\o_i$ and $\l$ are dominant we
conclude that $w \o_i = \o_i$ which is a contradiction. \qed

\

Now we will prove the main theorem.

\begin{thm10}
Let $F$ be a Steinberg fiber of $G_{\s}$ and $\bar{F}$ its closure
in $X_{\s}$. Then $$\bar{F}-F=\bigsqcup_{J \subset I} \bigsqcup_{
\substack{w \in W^{\s(J)} \\ \supp_{\s}(w)=I}} Z^w_{J, \s},$$
which also coincides with the set $Z(a_1,\dots,a_l) \cap (X_{{\rm sc},\s}
- G_{{\rm sc},\s})$ for all
$a_1,\dots,a_l$.
\end{thm10}

Proof.  Let $C$ be an irreducible component of
$Z(a_1, \cdots, a_l)$. Then by Krull's principal ideal
theorem, $\dim(C) \ge \dim(G_{\rm sc})-l$. By
\cite[8.20]{L2},
$$\dim(Z^w_{J, \s})=\dim(G)-l(w)-|I-J|<\dim(G_{\rm sc})-l,$$ for
$J \neq I$ and $w \in W^{\s(J)}$ with $\supp_{\s}(w)=I$. Thus
by Lemma 7.2, $$\dim \bigl(C \cap (X_{\rm sc, \s}-G_{\rm sc, \s})
\bigr)<\dim(G_{\rm sc})-l \leq \dim(C).$$ Hence $C \cap G_{\rm sc,
\s}$ is dense in $C$. But, by Corollary 6.3, the intersection
$Z(a_1, \cdots,a_l) \cap G_{\rm sc, \s}$ is a single
Steinberg fiber $F(a_1,\dots,a_n)$ which, as in 3.2, is irreducible. We
conclude that $C$ is contained in the
closure of $F(a_1, \cdots,a_l)$, and thus the closure of $F(a_1,
\cdots, a_l)$ is $Z(a_1, \cdots, a_l)$. In particular, $Z(a_1,
\cdots, a_l)$ is irreducible.

Let $F$ be a Steinberg fiber of $G_{\s}$. Then $F=\pi(F(a_1,
\cdots, a_l))$ for some $a_1, \cdots, a_l \in k$. Hence
$\bar{F}=\pi(Z(a_1, \cdots, a_l))$. The statement now follows from
Lemma 7.2 and Lemma 5.2. \qed

\begin{remark} 1. We call an element $w \in W$ a $\s$-twisted Coxeter
element if $l(w)=l$ and $\supp_{\s}(w)=I$. (The notation of
twisted Coxeter elements was first introduced by Springer in
\cite{Sp1}. Our definition is slightly different from his).
It follows easily from Theorem 7.3 that $\overline{Z^w_{I-\{i\},
\s}}$ are the irreducible components of $\bar{F}-F$, where $i \in
I$ and $w$ runs over all $\s$-twisted Coxeter elements that are
contained in $W^{I-\{\s(i)\}}$.

2. By the proof of Theorem 7.3 we may also deduce that the closure
of a Steinberg fiber $F$ within $X_{{\rm sc},\s}$ coincides with
$Z(a_1,\dots, a_l)$ for certain uniquely determined $a_1, \dots,
a_l$ depending on $F$. This result may be considered as an
extension of Corollary 2 in \cite{Sp3} to the compactification
$X_{{\rm sc,\s}}$ of $G_{\rm sc, \s}$. More precisely, notice that
the statement of \cite[Corollary 2]{Sp3} is equivalent to saying
that a Steinberg fiber $F$ of $ G_{\rm sc, \s}$ is the common zero
set of the functions $f_{i,a_i}^\s$ for uniquely determined
$a_1,\dots,a_l$. Here we think of $f_{i,a_i}^\s$ as regular
functions on $G_{\rm sc, \s}$ as explained in 6.1. When
generalizing to $X_{\rm sc,\s}$ the only difference is that we
have to regard $\bar{f}_{i,a_i}^\s$ as sections of certain line
bundles on $X_{\rm sc, \s}$.
\end{remark}

\

Similar to \cite[4.6]{H1}, we have the following consequence.

\begin{thm11} Assume that $G_{\s}$ is defined and split over $\mathbb
F_q$, then for any Steinberg fiber $F$ of $G_{\s}$, the number of
$\mathbb F_q$-rational points of $\bar{F}-F$ is
$$(\sum_{w \in W} q^{l(w)}) (\sum_{\supp_{\s}(w)=I} q^{l(w_0
w)+L(w_0 w)}),$$ where $w_0$ is the maximal element of $W$ and for
$w \in W$, $l(w)$ is its length and $L(w)$ is the number of simple
roots $\a$ satisfying $w \a<0$.
\end{thm11}

\section*{Acknowledgements}
We thank J.C.Jantzen, G.Lusztig and T.A.Springer for some useful
discussions and comments.

\bibliographystyle{amsalpha}

\end{document}